\documentclass{amsart}

\newtheorem{theorem}{Theorem}[section]
\newtheorem{lemma}[theorem]{Lemma}

\theoremstyle{definition}
\newtheorem{definition}[theorem]{Definition}
\newtheorem{example}[theorem]{Example}

\theoremstyle{remark}
\newtheorem{remark}[theorem]{Remark}

\numberwithin{equation}{section}

\begin{document}

\title{Generalized Sasakian Structures from a Poisson Geometry Viewpoint}

\author{Janet Talvacchia}
\address{Department of Mathematics and Statistics \\500 College Ave\\ Swarthmore, PA 19081}
\curraddr{}
\email{jtalvac1@swarthmore.edu}
\thanks{}




\date{May 2, 2023}


\begin{abstract} In this paper we define a canonical Poisson structure on a normal generalized contact metric space and use this structure to define a generalized Sasakian structure. We show also that this canonical Poisson structure enables us to distinguish generalized Sasakian structures from generalized coK\"ahler structures.
\end{abstract}

\maketitle

\section{Introduction}\label{S:intro} The notion of a generalized complex structure, introduced by Hitchin in his paper \cite{[H]} and developed by Gualtieri (\cite{[G1]},\cite{[G2]}) gives a framework that unifies both complex and symplectic structures on even dimensional manifolds. The odd dimensional analog of this structure, a generalized contact structure, was taken up by  Vaisman (\cite {[V1]},\cite{[V2]}), Poon, Wade \cite{[PW]}, Sekiya \cite{[S]}, and Aldi and Grandini \cite{[AG]}. This framework unifies almost contact, contact, and cosymplectic structures.  In odd dimensions, classical coK\"ahler and classical Sasakian structures are special almost contact metric structures distinguished by the fact that both generate K\"ahler structures on the cone $M\times \mathbb{R}$.  It is natural to consider what  their generalized counterparts would be.

\par
Generalized K\"ahler structures were introduced by Gualtieri (\cite{[G1]},\cite{[G2]},\cite{[G3]}) and found their way immediately into the physics literature (\cite{[Hu]}, \cite{[LMTZ]}).   Generalized coK\"ahler spaces were defined and studied by the author and R.Gomez in \cite{[GT2]} where it was shown that the standard product construction yielded a generalized K\"ahler manifold if and only if each of the factors in the product was a generalized coK\"ahler manifold.  It is shown there as well that classical Sasakian manifolds lie outside of this class of spaces.  (Generalized coK\"ahler spaces were later studied by Wright \cite{[W]} from the point of view of reduction of generalized K\"ahler spaces.) In \cite{[T]} the author showed that, indeed, any notion of generalized Sasakian could not arise via a standard product construction.  What a possible notion of generalized Sasakian might be has been discussed in (\cite{[IW]}, \cite{[S]}, \cite{[V1]}, \cite{[V2]}) but the definitions, while encapsulating many of the properties one would want in generalized Sasakian space,  are not all consistent, are defined on the cone $M \times \mathbb{R}$  instead of intrinsically on the space $M$,  and don't give an immediate way to distinguish generalized Sasakian from generalized coK\"ahler spaces.

\par
In  Abouzaid and Boyarchenko \cite{[AB]}, an approach to generalized geometry via Poisson structures was written down. They showed that there was a canonical Poisson structure associated to a generalized complex structure and this Poisson structure in fact determines the generalized complex structure.   (The existence of such a Poisson structure was independently noted  by Gualtieri \cite{[G3]} and Lyakhovich and Zabzine\cite{[LZ]}.)  In more recent years, several authors (see for example \cite{[C]}, \cite{[B]}, \cite{[BG]}, \cite{[G4]}) approach generalized complex geometry from the point of view of the underlying Poisson structure as the fundamental object.  An approach to generalized contact structures via Poisson geometry has not previously been undertaken.  The aim of this paper is two-fold.  The first is to define a canonical Poisson structure on a normal generalized contact structure that characterizes it.   We do that in section 3.  The second objective is to use this structure to propose a definition of a generalized Sasakian structure in terms of an invertibility criterion on the Poisson structure.  We do that in sections 4 and 5.   In \cite{[G4]}, Gualtieri proved results that describe how a generalized K\"ahler structure can be generated via gauge transformations of the Poisson structure underlying a generalized complex structure.  In section 4, we look at these theorems in the context of $M \times \mathbb{R}$ where $M$ is a normal generalized contact metric manifold whose underlying Poisson structure is the  product Poisson structure of canonical Poisson structure on $M$ defined earlier in the paper and the canonical Poisson structure on  the space $\mathbb{R}$  with  its usual normal generalized contact metric structure.  In section 5, we then reduce the conditions on the Poisson structure on $M \times \mathbb{R}$ to conditions on the canonical Poisson structure on $M$.  This suggests a definition of generalized Sasakian in terms of an invertibility criterion on this Poisson structure on $M$.  We show that classical Sasakian spaces are always generalized Sasakian. Generalized coK\"ahler structures are shown to never satisfy this criterion so it  distinguishes between the generalized Sasakian and generalized coK\"ahler cases.  In essence, we show that  while generalized coK\"ahler spaces generate generalized K\"ahler spaces via a product construction,  generalized Sasakian spaces generate generalized K\"ahler spaces via gauge transformations of its underlying canonical Poisson structure.

\section{Preliminaries}\label{B:Back}

\indent We begin with a very brief review of generalized geometric structures. Throughout this paper we let $M$ be a smooth manifold. Consider the big tangent bundle, $TM\oplus ~T^*M$.  We define a neutral metric on $TM\oplus T^*M$ by$$  \langle X + \alpha  , Y + \beta \rangle =  \frac{1}{2} (\beta (X) + \alpha (Y) )$$ and the Courant bracket by $$[[X+\alpha, Y+ \beta ]] = [X,Y] + {\mathcal L}_X\beta -{\mathcal L}_Y\alpha -\frac{1}{2} d(\iota_X\beta - \iota_Y\alpha)$$ where $X, Y \in TM$ and $\alpha ,\beta  \in T^*M$ . A subbundle of $TM\oplus T^*M$ is said to be involutive  or integrable if its sections are closed under the Courant bracket\cite{[G1]}.

\begin{definition} \cite{[G1]}
A generalized  almost complex structure on $M$ is an endomorphism $\mathcal J$ of $TM\oplus T^*M$ such that $\mathcal J + \mathcal J^* =  0 $ and $\mathcal J^2 = - Id$. If the $\sqrt{-1}$ eigenbundle $L\subset (TM\oplus TM^{*})\otimes\mathbb{C}$ associated to $\mathcal J$ is involutive with respect to the
Courant bracket, then $\mathcal J$ is called a generalized complex structure.
\end{definition}

\indent Here are the prototypical examples:

\begin{example}  \cite {[G1]}
Let $(M^{2n}, J)$ be a complex structure.  Then we get a generalized complex structure by setting
$$\mathcal J_{J} = \left ( \begin{array}{cc}  -J & 0 \\ 0 & J^* \end{array} \right ).$$
\end{example}

\begin{example}  \cite {[G1]}
Let $(M^{2n}, \omega )$ be a symplectic structure.  Then we get a generalized  complex structure by setting
$$\mathcal J_{\omega} = \left ( \begin{array}{cc}  0 & -\omega^{-1} \\ \omega & 0 \end{array} \right ).$$
\end{example}
\indent Diffeomorphisms of $M$ preserve the Lie bracket of smooth vector fields and in fact diffeomorphisms are the
only such automorphisms of the tangent bundle. But in generalized geometry, there is actually more flexibility. 
Given $TM\oplus T^{*}M$ equipped with the Courant bracket, the automorphism group is comprised of the diffeomorphisms of $M$
and some additional symmetries called \emph{B}-field transformations \cite{[G1]}.
\begin{definition}\cite{[G1]}
Let $B$ be a two-form which we view as a map from $TM \rightarrow T^{*}M$ given by interior product. Then the invertible bundle map
$$e^{B}:= \left ( \begin{array}{cc}  1 & 0 \\ B & 1 \end{array} \right):X+\xi \longmapsto X+\xi + \iota_{X}B$$
is called a B-field transformation.
\end{definition}
 The bundle map defined above by the two form $B$ preserves the Courant bracket if and only if $B$ is closed (see \cite{[G1]}).

\indent Recall that we can reduce the structure group of $TM\oplus T^{*}M$ from $O(2n,2n)$ to the maximal compact subgroup $O(2n)\times O(2n)$. This
is equivalent to an orthogonal splitting of $TM\oplus T^{*}M=V_{+}\oplus V_{-}$, where $V_{+}$ and
$V_{-}$ are positive and negative definite respectfully with respect to the inner product. Thus we can define a positive definite Riemannian metric
on the big tangent bundle by $$G=<,>|_{V_{+}}-<,>|_{V_{-}}.$$
\begin{definition}\cite{[G1]}
A  generalized metric $G$ on $M$ is an automorphism of $TM\oplus T^*M$ such that $G^{*}=G$ and $G^{2}=1.$ 
\end{definition}

In the presence
of a generalized almost complex structure $\mathcal J_1$, if $G$ commutes with $\mathcal J_1$ ($G\mathcal J_1 = \mathcal J_1  G$) then $G\mathcal J_1 $ squares to $-1$ and we generate a second generalized almost complex structure, $\mathcal J_2$ $= G\mathcal J_1$, such that $\mathcal J_1$ and $\mathcal J_2$ commute and $G=-\mathcal J_1 \mathcal J_2$.
This leads us to the following:
\begin{definition}\cite{[G1]}A  generalized K\"ahler structure is a pair of commuting generalized complex structures $\mathcal J_{1}, \mathcal J_{2}$ such that
$G=-\mathcal J_{1}\mathcal J_{2}$ is a positive definite metric on $T\oplus T^{*}.$
\end{definition}

\begin{example} \cite{[G1]}
Consider a K\"ahler structure $(\omega,J,g)$ on $M$. By defining $\mathcal J_{J}$ and $\mathcal J_{\omega}$ as in
Examples 2.2 and 2.3, we obtain a generalized K\"ahler structure on $M$, where
$$G=\left ( \begin{array}{cc}  0 & g^{-1} \\ g & 0 \end{array} \right ).$$

\end{example}
\indent We now recall the odd dimensional analog of generalized complex geometry. We use the definition given by Sekiya  (see \cite{[S]}).
\begin{definition} \cite{[S]} A generalized almost contact structure on $M$ is a triple $(\Phi, E_\pm)$ where $\Phi $ is an endomorphism of $TM\oplus T^*M$, and $E_+$ and $E_-$ are sections of $TM\oplus T^*M$ which satisfy
\begin{equation}
\Phi + \Phi^{*}=0
\end{equation}
\begin{equation}\label{phi}
\Phi \circ \Phi = -Id + E_+ \otimes E_- + E_- \otimes E_+
\end{equation}
\begin{equation}\label{sections}
 \langle E_\pm, E_\pm \rangle = 0,  \  \    2\langle E_+, E_-  \rangle = 1.
\end{equation}

\end{definition}
Now, since $\Phi$ satisfies $\Phi^3 + \Phi =0$, we see that $\Phi$ has $0$ as well as $\pm \sqrt{-1}$ eigenvalues when viewed as an endomorphism of the complexified big tangent bundle $(TM\oplus ~T^*M) ~\otimes { ~\mathbb C}$.  The kernel of $\Phi$ is $L_{E_+} \oplus L_{E_-}$ where $L_{E_\pm}$ is the line bundle spanned by ${E_\pm}$.  Let $E^{(1,0)}$ be the $\sqrt{-1}$ eigenbundle of $\Phi$.  Let $E^{(0,1)}$ be the $-\sqrt{-1}$ eigenbundle. Observe:

$$
E^{(1,0)} = \lbrace X + \alpha - \sqrt{-1}  \Phi (X + \alpha ) |  \langle E_\pm, X + \alpha \rangle = 0 \rbrace
$$

$$
E^{(0,1)} = \lbrace X + \alpha + \sqrt{-1}  \Phi (X + \alpha ) | \langle E_\pm, X + \alpha \rangle = 0 \rbrace .$$

Then the complex vector bundles
$$L^+ = L_{E_+} \oplus E^{(1,0)}$$
and
$$L^- = L_{E_-} \oplus E^{(1,0)}$$
are maximal isotropics.
\begin{definition} \cite{[PW]}
A generalized almost contact structure $(\Phi,E_{\pm})$ is a  generalized contact structure if either $L^{+}$ or $L^{-}$ is closed with respect to the Courant bracket. The generalized
contact structure is strong if both $L^{+}$ and $L^{-}$ are closed with respect to the Courant bracket.
\end{definition}

\begin{definition} \cite{[GT2]}
A  generalized almost contact structure $(M, \Phi,E_{\pm})$ is a normal generalized contact structure if  $\Phi$ is strong and $[[E_+, E_-] ] = 0$.
\end{definition}

\begin{remark}  This definition of normality is motivated by Theorem 1 of \cite{[GT1]} that shows that product of two generalized almost contact spaces  $(M_1, \Phi_1, E_{\pm ,1})$ and  $(M_2, \Phi_2, E_{\pm , 2})$  induces a standard generalized almost  complex structure on $M_1 \times M_2 $.  The generalized complex structure is integrable if each $\Phi_i$ is strong and   $[[E_{+,i},  E_{-,i}]] = 0$. 
\end{remark}
Here are the standard examples:

\begin{example}\label{almost contact example}  \cite{[PW]}
Let $(\phi , \xi, \eta)$ be a normal almost contact structure on a manifold $M^{2n+1}$.  Then we get a generalized almost contact structure by setting
$$ \Phi = \left ( \begin{array}{cc}  \phi & 0 \\ 0 & -\phi^* \end{array} \right ),\  \   E_+ = \xi,\  \   E_-= \eta $$  where $(\phi^*\alpha )(X) = \alpha (\phi (X)), \   X \in TM,\   \alpha \in T^{*}M$. Moreover, $(\Phi, E_\pm)$ is an example of a strong generalized contact structure.  In fact, this structure is a normal generalized contact structure.
\end{example}

\begin{example} \label{contact example} \cite{[PW]}
Let $( M^{2n+1}, \eta )$ be a contact manifold with $\xi $ the corresponding Reeb vector field so that
$$ \iota_\xi d\eta = 0 \  \  \  \eta ( \xi ) = 1.$$
Then $$\rho ( X) := \iota_X d\eta - \eta ( X)\eta$$ is an isomorphism from the tangent bundle to the cotangent bundle.  Define a bivector field by
$$\pi (\alpha , \beta ) := d\eta (\rho^{-1}( \alpha ), \rho^{-1}( \beta )),$$
where $\alpha, \beta \in T^{*}M. $
We obtain a generalized almost contact structure by setting
$$ \Phi = \left ( \begin{array}{cc}  0 & \pi \\ d\eta & 0 \end{array} \right ),\  \   E_+ = \eta,\  \   E_-= \xi .$$
In fact, $(\Phi, E_{\pm})$ is an example of a generalized contact structure which is not strong.
\end{example}

\begin{definition} \cite{[S]}

A generalized almost contact metric structure is a generalized almost contact structure
$(\Phi , E_{\pm})$ along with a generalized Riemannian metric $G$ that satisfies
\begin{equation}\label{compatcond}
-\Phi G \Phi = G - E_+ \otimes E_+ -E_- \otimes E_-.
\end{equation}

\end{definition}

\begin{definition} \cite{[GT2]}

A generalized coK\"ahler structure is a normal generalized contact metric structure  $(M, \Phi, E_+ , E_- , G)$ where both $\Phi$ and $G\Phi$ are strong.

\end{definition}

The following theorem was proved in \cite{ [GT2]}

\begin{theorem}\label{T1}
Let $M_{1}$ and $M_{2}$ be odd dimensional smooth manifolds each with a  generalized contact metric structure
$(\Phi,E_{\pm,i},G_i),i=1,2$ such that on the product $M_1\times M_2$ are two
 generalized almost complex structures: $\mathcal{J}_1$ which is the natural generalized almost complex structure induced from $\Phi_1$ and $\Phi_2$ and $\mathcal{J}_2=G\mathcal{J}_1$ where
$G=G_1\times G_2$. Then $(M_1\times M_2,\mathcal{J}_1,\mathcal{J}_2)$ is generalized K\"ahler if and
only if $(\Phi_{i},E_{\pm,i},G_i)$, $i=1,2$ are  generalized coK\"ahler structures.
\end{theorem}

\section{The Canonical Poisson Structure on a Generalized Almost Contact Manifold}\label{P:Poisson}

\indent In this section, we show there is a canonical Poisson structure underlying a strong generalized almost contact structure $(M, \Phi, E_+ , E_- )$ provided $[[E_+ , E_- ]] = 0$.  The condition $[[E_+ , E_- ]] = 0$ arose previously in \cite{[GT1]}  as a necessary and sufficient condition for the product of two strong generalized contact manifolds to be a generalized complex manifold.  It is not surprising that it appears as a condition to obtain a Poisson structure on the generalized almost contact space that is formed from the defining elements of the space.

\par

First, we show that  a manifold $M$ with generalized contact structure $(M, \Phi, E_+ , E_- )$ satisfying the condition $[[E_+ , E_-]]=0$ admits a foliation.  Let $e_\pm = pr_{TM}E_\pm$.  Since $[[E_+, E_-]]=0$, it follows that $[e_+, e_-] = 0$.  The Frobenius theorem implies that $\lbrace e_+, e_-\rbrace$ is an integrable distribution that induces a foliation of $M$.  If $(M^{2n+1}, \Phi, E_+ , E_- )$  is in fact a strong generalized contact structure, then  the foliation has one dimensional leaves and in a neighborhood of a point $p$, $(E_+ \oplus E_-)^\perp = E^{(1,0)}\oplus E^{(0,1)}$ has a generalized complex structure transverse to this foliation. We see this as follows.  Recall that a complex distribution $D \subset TM\otimes {\mathbb C} = T_{\mathbb C}M$  of constant complex codimension $k$ on a real manifold $M$  of dimension $n$ is integrable if, in some neighborhood $U$ of each point $m \in M$, there exists complex functions $f_1, \dots f_k \in C^\infty (U, {\mathbb C})$ such that $\lbrace df_1, \dots , df_k\rbrace$ are linearly independent at each point in $U$ and annihilate all complex vector fields lying in $D$.  By the Newlander-Nirenberg theorem \cite { [NN]} ( see also \cite {[G1]} for an exposition of this version of the theorem), we know that $D$ is integrable if  $D$ is involutive, dim $D \cap \bar D$ is constant, and $D\oplus \bar D$ is involutive.  If these conditions are satisfied, every point $m \in M$ has a neighborhood $U$ isomorphic, as a smooth manifold with complex distribution, to an open set in ${\mathbb R}^{n-2k} \times {\mathbb C}^k$ which has  a natural distribution spanned by $\lbrace \frac {\partial}{\partial x_1}, \dots , \frac {\partial}{ \partial x_{n-2k}}, \frac {\partial}{\partial z_1}, \dots , \frac {\partial } {\partial z_k}\rbrace $. 

In the situation we are considering, we take $D = pr_{T_{\mathbb C}M}(E^{(1,0)})$.  
Since $(M^{2n+1}, \Phi, E_+ , E_- )$ is a strong generalized contact structure,  we have that $[[E^{(1,0)}, E^{(1,0)}]] \subset E^{(1,0)}$ (see \cite{[GT2]}) which implies that $D$ is involutive on $T_{\mathbb C}M$ with respect to the Lie bracket.  $D\cap {\bar D} = 0$  so dim $D \cap \bar D$ is constant.  Lastly,  $[[E^{(0,1)}, E^{(0,1)}]] \subset E^{(0,1)}$  as well so that $D\oplus \bar D$ is involutive and all the conditions of the Newland Nirenberg theorem are satisfied.  Let $U$ be the neighborhood of $p\in M$ isomorphic to ${\mathbb R} \times {\mathbb C}^n$ given by the theorem.  Let $W = pr_{\mathbb C^n} (U)$.  Note that $T_{\mathbb C}W\oplus T_{\mathbb C}^*W \approx  E^{(1,0)}\oplus E^{(0,1)}$ so that $\Phi$ restricted to $T_{\mathbb C}W\oplus T_{\mathbb C}^*W$ satisfies $\Phi^2 = -  id$.  That $\Phi$ is integrable follows from its properties as part of a strong generalized contact structure.  Thus  the neighborhood $W$ has a canonical Poisson structure $\pi_0$ by the results of Abouzaid and Boyarchenko \cite {[AB]}. We can extend this bivector $\pi_0$ on $W$ to a bivector $\pi_M$  on $U \subset M$  as follows.  Let $e_+ = pr_{TM}( E_+)$ and let $e_- = pr_{TM}( E_-)$.  At each $ p \in M$, define $$\pi_M = \pi_0 + e_+ \wedge e_- .$$

\begin{lemma}
If  $[[E_+ , E_- ]] = 0$ then $\pi_M$ is a Poisson structure on $M$ .  

\end{lemma}
\begin{proof}
Note  that $[[E_+ , E_- ]] =0$ implies $\lbrack e_+ , e_-\rbrack = 0$.

 Now 
$$\lbrack \pi_M , \pi_M \rbrack  =  \lbrack \pi_0 , \pi_0 \rbrack + 2 \lbrack  \pi_0, e_+ \wedge e_- \rbrack  +\lbrack e_+ \wedge e_- , e_+ \wedge e_- \rbrack .$$

$ \lbrack \pi_0 , \pi_0 \rbrack = 0$ since $\pi_0 $ is a Poisson structure.

$  \lbrack  \pi_0, e_+ \wedge e_- \rbrack = \lbrack \pi_o , e_+ \rbrack \wedge e_ +  + e_+ \wedge \lbrack \pi_0 , e_-\rbrack $
We see that $\lbrack \pi_o , e_\pm \rbrack = 0$ by using the definition of the Schouten bracket and the fact that $\pi_0$ is a bivector on $(E_+ \oplus E_-)^\perp$. Thus,  $\lbrack  \pi_0, e_+ \wedge e_- \rbrack = 0$

The condition   $\lbrack e_+ , e_-\rbrack = 0$ insures  $\lbrack e_+ \wedge e_- , e_+ \wedge e_- \rbrack = 2 \lbrack e_+ ,e_-\rbrack \wedge e_+ \wedge e_-  = 0.$   

\end{proof}

\begin{example} 
  Let  $ M=SU(2)$.
On the Lie algebra $su(2)$ choose a basis $\lbrace X_1,  X_2,  X_ 3 \rbrace $ and a dual basis $\lbrace  \sigma^1, \sigma^2, \sigma^  3 \rbrace$ such that $[ X_i, X_j ] = -X_k$ and $ d\sigma^i = \sigma^j \wedge \sigma^k$ for cyclic permutations of $\lbrace i, j, k \rbrace$.  One can construct a classical normal almost contact structure by taking $\phi = X_2 \otimes \sigma^1 - X_1 \otimes \sigma^2$, $ \xi = X_3$, and $\eta = \sigma^3$.  Then, as in Example 2.12, we can construct a generalized almost contact structure by letting
$$ \Phi = \left ( \begin{array}{cc}  \phi & 0 \\ 0 & -\phi^* \end{array} \right ), \    E_+ = X_{3},\    E_-= \sigma^{3} $$  where $(\phi^*\alpha )(X) = \alpha (\phi (X)), \   X \in TM,\   \alpha \in T^{*}M$. 

 We compute:
 
 $$E^{(1,0)}_\phi = {\rm span} \lbrace X_1- \sqrt{-1}X_2 ,  \sigma^1- \sqrt{-1} \sigma^2 \rbrace$$
 $$E^{(0,1)}_\phi = {\rm span} \lbrace X_1+ \sqrt{-1}X_2 ,  \sigma^1 + \sqrt{-1} \sigma^2 \rbrace$$
 $$L^+ = {\rm span} \lbrace X_3, X_1 - \sqrt{-1}X_2 ,  \sigma^1- \sqrt{-1} \sigma^2 \rbrace$$ 
 $$L^{-} =  {\rm span} \lbrace \sigma^3, X_1 - \sqrt{-1}X_2 ,  \sigma^1 -\sqrt{-1} \sigma^2 \rbrace$$
Computing the Courant brackets we get

$$[[X_1 \pm \sqrt{-1}X_2,\sigma^1 \pm \sqrt{-1} \sigma^2]]=0$$
$$ [[X_{3},\sigma^{1}-\sqrt{-1}\sigma^{2}]]=-\sqrt{-1}(\sigma^{1}-\sqrt{-1}\sigma^{2})$$
$$[[X_{3}, X_{1}-\sqrt{-1}X_{2}]]=\sqrt{-1}(X_{1}-\sqrt{-1}X_2)$$
$$[[\sigma^{3},X_{1}-\sqrt{-1}X_2]]=\sqrt{-1}(\sigma^{1}-\sqrt{-1}\sigma^{2})$$ 
$$[[\sigma^{3},\sigma^1-\sqrt{-1}\sigma^2]]=0$$ 
\medskip

Also we compute that $[[E_{+},E_{-}]]=\mathcal{L}_{X_3}\sigma^{3}=0.$ Thus we see that $L_+$ and $L_-$ are closed under the Courant bracket, $\Phi_\phi $ is strong and  $ (\Phi_\phi, E_{\pm})$ is a normal generalized contact structure on $SU(2)$. One can verify directly that $[[E^{(1,0)}, E^{(1,0)}]] \subset E^{(1,0)}$ and $[[E^{(0,1)}, E^{(0,1)}]] \subset E^{(0,1)}$.  Since  $e_+ = X_3$ and $e_- = 0$, the foliation of $M$ is the characteristic foliation induced by the Reeb vector field $X_3$.  Any point $p \in M$ has local coordinates $x_1, z_1$  in some neighborhood $U$ where $\frac {\partial} {\partial x_1 }= X_3$ and  $\frac {\partial} {\partial z_1 } = X_1 - \sqrt{-1} X_2$.  $\Phi_\phi$ restricted to $W= pr_{\mathbb C} (U)$ defines a generalized complex structure which has a canonical Poisson structure $\pi_0$. At any point $p\in SU(2)$ we define $\pi_M = \pi_0 + e_+ \wedge e_ -$. Note that $e_+ \wedge e_ - = 0$ here as one can observe happens with the generalized contact structures arising from classical structures.

\end{example}

\begin{remark}
  Vaisman has constructed important Poisson structures associated to generalized CRF structures  (see \cite {[V2], [V3], [V4]}).  When $(M, \phi, E_+,  E_-)$ is a quasi-classical generalized CRF structure that admits a foliation, the Poisson structure here coincides with Vaisman's.  (See  \cite {[V4]}, section 3, for a discussion of quasi-classical generalized CRF structures with foliations.) It would be interesting to study further possible relationships between these Poisson structures.   

\end{remark}

 \par
In the rest of this paper, we look at this canonical Poisson structure on normal generalized contact metric spaces.

\section{Gauge Transformations of Poisson structures and Gualtieri's Theorem}\label{G: Gauge}
\indent In this section we review a theorem by Gualtieri \cite{[G4]} that we will use to generate the invertibility criterion on the canonical Poisson structure of the normal generalized contact metric manifold.  Recall that a Dirac structure is a maximal isotropic subbundle  $L \subset  TM \oplus T^*M$ that is involutive under the Courant bracket. Let $L_1$ and $L_2$ be transverse Dirac structures, $pr_{TM}(L_1) + pr_{TM}(L_2)= TM$.  Then their sum is defined by
$$L_1 + L_2 = \lbrace X + \alpha + \beta \  | \   X+ \alpha \in L_1, X+\beta \in L_2\rbrace .$$
We use the notation $L_2-L_1$ for  the Dirac structure $(-1) L_1 + L_2$.

\noindent  The scaling of a Dirac structure by a nonzero number $\lambda \in \mathbb{R}$ is defined by
 
 $$ \lambda L = \lbrace X + \lambda \alpha \  | \   X + \alpha \in L  \rbrace. $$

\par
 A real Poisson structure $\pi \in C^\infty (\bigwedge^2TM)$ on a real smooth manifold $M$ can be viewed as a Dirac structure by considering its graph sub-bundle $\Gamma_\pi \subset TM\oplus T^*M$.
$$\Gamma_\pi = \lbrace \pi \xi + \xi \   | \   \xi \in T^*M\rbrace$$
Since $[\pi , \pi]= 0$, $\Gamma_\pi$ is an involutive maximal isotropic on $TM\oplus T^*M$.

\begin{definition}
Two real Poisson structures $\pi_0$ and $\pi_1$ on $M$  are gauge equivalent if there exists a real closed two form $B$ such that 
$$e^B\Gamma_{\pi_0} = \Gamma_{\pi_1}$$
\end{definition}

In \cite{[SW]}, where this notion was introduced, it was shown that $\pi_0$ and $\pi_1$ are gauge equivalent if and only if $(I + B\pi_0)$ is invertible as a bundle automorphism of $T^*M$.  When this happens, 
$$\pi_1 = \pi_0 (I + B\pi_0)^{-1}$$

We see from this that $\pi_0$ and $\pi_1$ are gauge equivalent when their symplectic leaves coincide.

\begin{theorem} (Gualtieri,  \cite{[G4]})

Fix a generalized complex structure $\mathcal{J}_0$ on an even dimensional manifold $M$ with $\sqrt{-1}$ eigenbundle $L_0$ and underlying real Poisson structure $\pi_0$.  Let $\beta$ be a complex two form.  Then $L_1 = e^\beta L_0$ defines a generalized complex structure $\mathcal{J }_1$ if and only if $(I + B\pi_0)$ is invertible for $B = {\rm Im} \ \beta$.  The underlying real Poisson structure $\pi_1$ of $\mathcal{J}_1$ is given by $\pi_1 = \pi_0(I + B\pi_0)^{-1}$.

\end{theorem}

Let's look at this theorem now in the context of an even dimensional space $M\times \mathbb{R}$ where $M$ is a normal generalized contact metric manifold and $\mathbb{R}$ is given its standard normal generalized contact metric structure $$( {\mathbb R}, \Phi_\mathbb{R} \equiv 0, (0, dt), (\frac {\partial }{ \partial t}, 0), \  G_\mathbb{R} =  \left ( \begin{array}{cc}  0 & (dt^2)^{-1} \\ dt^2 & 0 \end{array} \right )\  )$$ \    where $dt^2$ is the standard euclidean metric.

The canonical Poisson structure $\pi_\mathbb{R}$ underlying this generalized contact  metric structure on $\mathbb{R}$ is the zero Poisson structure.  We get a Poisson structure on $M\times \mathbb{R}$ by taking the product Poisson structure : $\pi = \pi_M  + \pi_\mathbb{R}$.  Theorem 1 of \cite{[GT1]} implies that $M\times \mathbb{R}$ admits a generalized complex structure ${\mathcal J}_0$ constructed from the data on $M$ and $\mathbb{R}$.  In fact , an explicit formula for ${\mathcal J}_0$ is given there: 
\begin{align*}
{\mathcal J}_0 ( X+\alpha , a \frac{\partial}{ \partial t} + b dt) =  &( \Phi(X+\alpha) - 2\langle dt , a \frac{\partial}{ \partial t} + b dt \rangle E_+ - 2\langle  \frac{\partial }{\partial t} , a \frac{\partial}{ \partial t} + b dt\rangle E_- , \cr
 &0  + 2\langle E_+, X +\alpha  \rangle dt + 2\langle E_- ,  X +\alpha \rangle \frac{\partial}{ \partial t} ).
 \end{align*}
 
 From this formula we see that $\pi$ is the underlying Poisson structure for $\mathcal{J}_0$.  Let $\eta = pr_{T^*M}(E_+) + pr_{T^*M}(E_-)$.  Assume $d\eta \neq 0$.  Let $B = d(e^t\eta)$.  Then if $(I + B\pi)$ is invertible, we get a new complex structure $\mathcal{J}_1$ on $M\times \mathbb{R}$.  To show that ${\mathcal{J}_0}$ and ${\mathcal J}_1$ form a generalized K\"ahler structure, we must show that ${\mathcal J}_0$ and $\mathcal{J}_1$ commute and $\mathbb{G} = -\mathcal{J}_0$$\mathcal{J}_1$ is positive definite.  In \cite{[G4]}, Gualtieri proved  the following theorem that allows one to verify the generalized K\"ahler conditions on $\mathcal{J}_0$ and $\mathcal{J}_1$ in terms of their underlying Dirac structures $L_0$ and $L_1$:

\begin{theorem}  (Gualtieri \cite{[G4]})

The pair of complex Dirac structures $(L_0, L_1)$ defines a generalized K\"ahler structure if and only if it satisfies all of the following conditions:
\begin{enumerate}
\item $L_0$ is transverse to $\bar{L}_0$; i.e., $pr_{T_{\mathbb{C}}M}(L_0) + pr_{T_{\mathbb{C}}M}(\bar{L}_0) = T_{\mathbb{C}}M $. Similarly  $L_1$ is transverse to $\bar{L}_1$.

\item The real Dirac structures
$$\Gamma_{\pi_1} = \frac {1}{ 2\sqrt{-1}}(L_0  - \bar{L}_0 ),\   \   \   \   \Gamma_{\pi_2} =  \frac {1}{2 \sqrt{-1}}(L_1 - \bar{L}_1),$$
define real Poisson structures, i.e. $\Gamma_{\pi_1}\cap TM = \Gamma_{\pi_2} \cap T M= 0$.

\item The complex Dirac structures
$$L_{{\sigma_+ }} =  \frac{1}{ 2\sqrt{-1}}(L_0 - L_1), \  \  \  \   L_{\sigma_-} = \frac {1}{2\sqrt{-1}} (L_0 - \bar{L}_1)$$
define holomorphic Poisson structures $(I_+, \sigma_+),  (I_-, \sigma_-)$ respectively, i.e. $T_\mathbb{C} M  = (L_{{\sigma_+ }}\cap T_\mathbb{C} M ) \oplus (\bar{L}_{{\sigma_+ }} \cap T_\mathbb{C} M) $ and $T_\mathbb{C} M = (L_{{\sigma_- }} \cap T_\mathbb{C} M ) \oplus (\bar{L}_{{\sigma_- }} \cap T_\mathbb{C} M) $.  

\item For all nonzero $u\in L_0\cap L_1$ we have  $\langle u,\bar{u} \rangle   > \  0$.

\end{enumerate}

\end{theorem}

We use this theorem to verify the generalized K\"ahler conditions since we have the existence of $\mathcal{J}_1$ but no formula for it to allow us to check its commutativity with $\mathcal{J}_0$ directly.  In \cite{[GT1]}, there is an explicit formula for $\mathcal{J}_0$  as well as  for generators of $L_0$.  The relation $ L_1 = e^\beta L_0$ allows us thus to compute generators for $L_1$.

From  \cite{[GT1]}, $L_0$ is generated by
\begin{itemize}
\item $(E^{(1,0)}, 0)$
\item $(E_{-1}, - \sqrt{-1}\frac {\partial}{ \partial t} )$
\item $(E_{+1}, -\sqrt{-1}dt)$
\end{itemize}

Setting  $\beta = \sqrt{-1}\   d(e^t\eta)$, we see that $L_1 = e^\beta  L_0$ is generated by
\begin{itemize}
\item $(E^{(1,0)},0) + \iota_X d(e^t \eta )$ where $X+\alpha \in E^{(1,0)}$
\item $(E_{-1}, - \sqrt{-1} \frac{\partial}{ \partial t} ) + \iota_{(e_+ - \sqrt{-1}\frac {\partial}{ \partial t}) } d(e^t \eta )$
\item $(E_{-1}, - \sqrt{-1} dt) + \iota_{e_-}d(e^t \eta )$
\end{itemize}

It is a straightforward computation to verify the conditions of Theorem 4.3  using the generators of $L_0$ and $L_1$.  As an example, we verify condition (2): 
Consider $\Gamma_{\pi_1} = \frac{1}{2\sqrt{-1}} (L_0 - \bar{L}_0)$.  Recall that $L_0- \bar{L}_0 = \lbrace X + \alpha - \beta\  | \  x+ \alpha \in L_0 , \  \  X + \beta \in \bar{L}_0 \rbrace$. In particular, note that the vector field component in $L_0- \bar{L}_0 $ comes from elements of $L_0$ and $\bar{L}_0$ with the same vector field part.  Consider, now, the generators of $L_0$ and $\bar{L}_0$.  Since $E^{(1,0)} \cap E^{(0,1)} = 0$ and $E_{\pm, 1} \in (E^{(1,0)} \oplus E^{(0,1)})^\perp$, the generator $(E^{(1,0)}, 0)$ contributes nothing to $\Gamma_{\pi_1}$.  The generator $(E_{-1}, - \sqrt{-1}\frac {\partial}{\partial t} ) \in L_0$  and its corresponding generator  $(E_{-1}, + \sqrt{-1}\frac {\partial}{ \partial t} ) \in \bar{L}_0$  also yield no terms in $\Gamma_{\pi_1}$.  The generators $(E_{+1}, -\sqrt{-1}dt)$ and 
$(E_{+1}, +\sqrt{-1}dt)$  contribute to $\Gamma_{\pi_1} $ so we see that $\Gamma_{\pi_1}  = \frac{1}{2\sqrt{-1}} {\rm span} \lbrace (E_{+,1},-2\sqrt{-1} dt )\rbrace$ and $\Gamma_{\pi_1}  \cap T_{\mathbb{C}} ( M\times {\mathbb R}) =\emptyset$.  The argument that $\Gamma_{\pi_2}  \cap T_{\mathbb{C}} ( M\times {\mathbb R}) =\emptyset$ goes similarly.

In summary, if $(I+ B\pi)$ is invertible for $B=d(e^t\eta)$  with $d\eta \neq 0$ and $ \pi = \pi_M + \pi_\mathbb{R}$,  we can use Gualtieri's theorems to see that a  normal generalized  contact metric structure generates a generalized K\"ahler structure on  $M\times \mathbb{R}$.

\section{ Reducing the Invertibility Condition to a Condition on $M$ and the Definition of Generalized Sasakian}\label {(M): Main}

Ideally, we'd like a condition directly on $M$ that would tell us whether this construction  of K\"ahler structures via a gauge transformation of the canonical Poisson structure on $M\times \mathbb{R}$ by $B = d(e^t \eta )$ will work.  We do this in the following theorem.

\begin{theorem}
Let $(M, \Phi , E_+ , E_- , G)$ be a normal generalized contact metric structure. Let $\pi_M$ be the canonical Poisson structure on $M$.  Let $\eta = pr_{T^*M}E_+ + pr_{T^*M}E_-$.

Let $(\   \mathbb{R}, \Phi = 0 , dt, \frac {\partial }{\partial t}, G_\mathbb{R} =  \left ( \begin{array}{cc}  0 & (dt^2)^{-1} \\ dt^2 & 0 \end{array} \right ) \  )$ be the standard normal generalized contact metric structure on $\mathbb{R}$.  Let $\pi_\mathbb{R}$ be  its underlying Poisson structure , $\pi_\mathbb{R} = 0$.

Let $\pi = \pi_M + \pi_\mathbb{R}$ be the product Poisson structure on $M\times \mathbb{R}$.

Let $B = d( e^t \eta)$. Let $\widehat{B}= d\eta $.

Then $(I+ B \pi)$ is invertible as a map from  $T^*( M \times \mathbb{R} )$ to $T^*( M \times \mathbb{R} )$ if and only if $(I + e^t \widehat B \pi_M)$ is invertible as a map from $T^*M$ to $T^*M$ for all $t\in\mathbb{R}$.

\end{theorem}
\begin{proof}

\begin{align*}
(I + B\pi ) &= I + B( \pi_M + \pi_\mathbb{R})\\
                &= I+B\pi_M  \\
                &= I + d(e^t \eta )( \pi_0 +\pi_1)\\
                &= I +e^t d\eta( \pi_0 +\pi_1 )+ e^tdt\wedge d\eta( \pi_0 +\pi_1)\\
\end{align*}
 Direct computation shows that $e^t dt \wedge d\eta (\pi_0 +\pi_1) =0$.   (To see this, write out $e^t dt \wedge d\eta (\pi_0 +\pi_1 )$as a sum of tensor products  using the definition of the wedge product and note that each summand  zero due to the absence  of any $\frac {\partial}{\partial t} $ terms.) So, $(I + B\pi ) =  I +e^t d\eta( \pi_0 +\pi_1 )$ as a map from  $T^*( M \times \mathbb{R} )$ to $T^*( M \times \mathbb{R} )$.
 
 Let $\lbrace x_1, x_2, \dots x_n, t\rbrace$ be local coordinates on $M\times \mathbb{R}$.  Let  $\alpha$ be an arbitrary 1-form in $T^*( M \times \mathbb{R} )$ given in local coordinates by $\alpha = a_1dx_1 + \dots a_mdx_n + a_tdt$.
 Then 
 \begin{align*}
 (I + B\pi )(\alpha ) &= (I + e^td\eta  (\pi_0+ \pi_1) )(a_1dx_1 + \dots  a_n dx_m +a_tdt) \\
 			     &= a_1dx_1 + \dots  a_n dx_n +a_tdt + e^td\eta  (\pi_0+ \pi_1) (a_1dx_1 + \dots  a_n dx_n +a_tdt) \\
			     &= a_1dx_1 + \dots  a_n dx_n +a_tdt + e^td\eta  (\pi_0+ \pi_1) (a_1dx_1 + \dots  a_n dx_n )  \\
			     &= (I+ e^t \widehat{B} \pi_M) (a_1dx_1 + \dots  a_n dx_n) + a_tdt
 \end{align*}
From this equation it is clear that $(I+ B \pi)$ is invertible as a map from  $T^*( M \times \mathbb{R} )$ to $T^*( M \times \mathbb{R} )$ if and only if $(I + e^t \widehat B \pi_M)$ is invertible as a map from $T^*M$ to $T^*M$ for all $t\in\mathbb{R}$.
\end{proof}

\begin{definition}
A normal generalized contact metric space $(M, \Phi , E_+ , E_- , G)$ is defined to be generalized Sasakian if $ (I + e^t(d\eta) \pi_M)$ is invertible as a map from $T^*M$ to $T^*M$  for all values of
 $t \in  \mathbb{R}$ where  $\pi_M$ is the canonical Poisson structure on $M$, $\eta = pr_{T^*M}E_+ + pr_{T^*M}E_-$ and $d\eta \neq 0$.
\end{definition}

 A natural question to ask is whether a  closed $B$ transform takes a generalized Sasakian structure to a generalized Sasakian structure.  Given a  generalized Sasakian structure $(M, \Phi, E_+, E_-, G)$, a closed $B$ transform takes it to a generalized contact metric structure $$(e^B\Phi e^{-B} , e^BE_+, e^BE_-, e^BGe^{-B})$$ satisfying $[[ e^B E_+, e^B E_-]] = 0$ (see \cite{ [S], [GT2]}). Let $E_\pm = X_\pm  + \alpha_\pm$. The following lemma shows that if $L_{X_+}B = 0$ and $L_{X_-}B = 0$, this $B$ transformed Sasakian structure is also generalized Sasakian.

\begin{lemma}
Let $(M, \Phi, E_+, E_-, G)$ be a generalized Sasakian structure. Let $E_+ = X_+ + \alpha_+$ and $E_- = X_- + \alpha_-$.    Let $B$ be a closed two form such that $L_{X_+}B = 0$ and $L_{X_-}B = 0$.  Then $(e^B\Phi e^{-B} , e^BE_+, e^BE_-, e^BGe^{-B})$ is a generalized Sasakian structure.
\end{lemma}
\begin{proof}  Since $(M, \Phi, E_+, E_-, G)$ is generalized Sasakian,  $ (I + e^t(d\eta) \pi_M)$ is invertible as a map from $T^*M$ to $T^*M$  for all values of
 $t \in  \mathbb{R}$ where  $\pi_M$ is the canonical Poisson structure on $M$, $\eta = pr_{T^*M}E_+ + pr_{T^*M}E_-$ and $d\eta \neq 0$.  Let $\eta_B = pr_{T^*M(}e^BE_+ )+ pr_{T^*M}(e^bE_-) = pr_{T^*M}E_+ + \iota_{X_+ }B + pr_{T^*M}E_- +\iota_{X_-}B$. Then $d\eta_B = d\eta +d \iota_{X_+ }B  + d \iota_{X_-}B = d\eta + L_{X_+}B +L_{X_-}B$.  Since $L_{X_+}B = 0$ and $L_{X_-}B = 0$, $d \eta = d \eta_B$.  Thus  $ (I + e^t(d\eta_B) \pi_M) =  (I + e^t(d\eta) \pi_M)$ is invertible as a map from $T^*M$ to $T^*M$  for all values of
 $t \in  \mathbb{R}$ as well and $(e^b\Phi e^{-B} , e^BE_+, e^BE_-, e^BGe^{-B})$ is Sasakian.

\end{proof}

\begin{theorem} A classical Sasakian space is generalized Sasakian.

\end{theorem}
\begin{proof}
Let $ M_1 = ( M, \phi , \eta , \xi, g)$ be a classical Sasakian space.  (Here, $\phi$ is a (1,1) tensor, $\xi$ is the Reeb vector field, $\eta$ is a 1-form, and $g$ is the metric.)

Let 
\begin{align*}
\Phi_1 =  \Phi_\phi &= \left ( \begin{array}{cc}  \phi & 0 \\ 0 & -\phi^* \end{array} \right ) \\    E_{+ 1}&= \xi \\  E_{-1}&= \eta 
\end{align*}
  be the generalized contact structure generated  by the classical almost contact structure $\phi$. The compatible generalized metric is $$G_1 =\left ( \begin{array}{cc}  0 & g^{-1} \\ g & 0 \end{array} \right ).$$
Since $e_+\wedge e_ - = 0$ for this space, $\pi_m= \pi_0$.  Also  for a classical Sasakian space, $\pi_0 $ restricted to  $(E_+ \oplus E_-)^\perp $ equals $(d \eta )^{-1}$.    Let $ \alpha $ be an arbitrary element of $T^*M$.   Then
 
 \begin{align*}
 (I+ e^t \widehat{B} \pi_M)(\alpha ) &= (I + e^t d\eta \circ \pi_0) (\alpha ) \\
                                      &= (1+e^t)(\alpha )  \\
 \end{align*}
 which is invertible.
\end{proof}

\begin{remark}   Definitions of generalized Sasakian structures have been put forth by Vaisman (\cite{[V1], [V2]}) and Sekiya (\cite{[S]}).  The definitions of Vaisman and Sekiya  are not invariant under B transforms hence this definition is distinct from theirs.

\end{remark}

Recall that for a classical coK\"ahler space $(M, \phi , \xi, \eta , g)$ we always have that $d\eta =0$.  So when we look at a classical coK\"ahler space as a generalized contact metric structure, the  gauge transform of its canonical poisson structure on $M \times \mathbb{R}$ by $B = d(e^t\eta )$ does not yield a new generalized complex structure. The map $( I + B \pi )$ reduces to just the identity map.  As a result we see that the  procedure using gauge transformation by $d(e^t\eta )$ on the canonical Poisson structure fails to generate a K\"ahler structure on $M\times \mathbb{R}$ even though a K\"ahler structure can be constructed on $M \times \mathbb{R}$ via a product construction.    This holds true with generalized coK\"ahler spaces as well.

\begin {theorem} Let  $(M, \Phi , E_+, E_- , G)$ be a generalized coK\"ahler space.  Let $\eta  = pr_{T^*M }(E_+)  + pr_{T^*M }(E_-)$.  Then $d\eta  = 0 $  when restricted to $ker  \  \eta = (E_+ \oplus E_-)^\perp$.   As a consequence,  we see that there is no intersection between the set of generalized Sasakian spaces and the set of generalized coK\"ahler spaces.

\end{theorem}

\begin{proof}  

Note that since $G(E_\pm) = E_\mp$,  the 1-form $\eta $ is the same for both generalized contact metric structures $( M, \Phi , E_+, E_-, G)$ and $( M, \tilde{\Phi} = G\Phi ,  \tilde{E}_+ = GE_+, \tilde{E}_- =GE_-, G)$. Since $\Phi$ and $G\Phi$ are both strong, the pair $(L , L^* )$ forms a Lie bialgebroid where $L = L_{E+} + E^{(1,0)}$ and $L^* =  L_{E-} + E^{(0,1)}$.  By theorem 2.9 of \cite{[PW]}, this is true if and only if  $d\eta  = 0 $  when restricted to $ker  \  \eta = (E_+ \oplus E_-)^\perp$.

\end{proof}

\section{ Arbitrary Products of Generalized Sasakian Manifolds}

The famous result of Calabi and Eckman  (see \cite{[CE]}) showed that the product of Sasakian spaces need not be K\"ahler.  This holds true for generalized Sasakian spaces as well.  One expects this since $(I + B(\pi_1 +\pi_2))$ need not be invertible even if both $ (I + B(\pi_1))$ and $ (I + B(\pi_2))$ are both invertible.  The details of the argument in section 4 shows that it works in part because the underlying Poisson structure on $\mathbb{R}$ is the zero Poisson structure.   The argument on the product  of generalized Sasakian spaces being K\"ahler will go through if one of the generalized contact metric structures has the zero Poisson structure as its canonical Poisson structure.  Hence, we get the following theorem:

\begin{theorem}
Let  $(M_1, \Phi , E_{+1} , E_{- 1}, G_1)$ be a generalized Sasakian space with canonical Poisson structure $\pi_{M_1}$.  
Let  $(M_2, \Phi , E_{+2} , E_{- 2}, G_2)$ be a generalized Sasakian space with canonical Poisson structure $\pi_{M_2 } = 0$.
Then $M_1 \times M_2$  admits a generalized K\"ahler structure.

\end{theorem}

Disclosures:  

 1. Conflict of interest :  There is no conflict of interest to disclose.
 
 2. Funding : There are no external funding sources to disclose.
 
3. Authors Contribution: Janet Talvacchia is the sole author of this paper. 

4. Data Availability statement :  There are no data used in this paper.

\end{document}